\documentclass [11pt,a4paper]{amsart}
\usepackage{amsmath}
\usepackage{stmaryrd}
\usepackage{textcomp}
\usepackage{amsmath}
\usepackage{cases}
\usepackage{amscd}
\usepackage{epsfig}
\usepackage{graphicx}

\usepackage{color}

\usepackage{verbatim}

\usepackage{amssymb}
\usepackage{amsfonts}
\usepackage{amsbsy}

\usepackage{amsthm}
\usepackage{amstext}
\usepackage{amsopn}
\usepackage{hyperref}
\hypersetup{colorlinks=true,
 linkcolor=blue,
filecolor=magenta,urlcolor=blue}

\newcommand{\calH}{\mathcal{H}}

\newcommand{\calC}{\mathcal{C}}

\newcommand{\calU}{\mathcal{U}}

\newcommand{\C}{\mathbb{C}}
\newcommand{\R}{\mathbb{R}}

\newcommand{\J}{\mathbb{J}}

\newcommand{\frakH}{\mathfrak{H}}

\newcommand{\bH}{{\bf H}}
\newtheorem{thm}{Theorem}[section]

\newtheorem{cor}[thm]{Corollary}

\newtheorem{prop}[thm]{Proposition}

\usepackage{fancyhdr}
\begin{document}
\title[Geodesics on the K\"ahler cone of the Heisenberg group]{Geodesics on the K\"ahler cone of the Heisenberg group}
\author{Joonhyung Kim \& Ioannis D. Platis \& Li-Jie Sun}

\address{Department of Mathematics Education,
Chungnam National University,
 99 Daehak-ro, Yuseong-gu, Daejeon 34134,
 Korea.}
\email{calvary@cnu.ac.kr}

\address{Department of Mathematics and Applied Mathematics,
University of Crete,
 Heraklion Crete 70013,
 Greece.}
\email{jplatis@math.uoc.gr}
\address {Department of applied science, Yamaguchi University,
2-16-1 Tokiwadai, Ube 7558611,
Japan.
}
\email{ljsun@yamaguchi-u.ac.jp}

\keywords{Heisenberg group, K\"ahler cone, geodesics.\\
\;\;{\it 2020 Mathematics Subject Classification:} 53C22, 32Q15.}

\thanks{{\it Acknowledgments.}  Parts of this work have been carried out while IDP was visiting Hunan University, Changsha, PRC and JK was visiting University of Crete, Greece and Yamaguchi University, Japan. Hospitality is gratefully appreciated. JK was supported by the National Research Foundation of Korea (NRF) grant funded by the Korea government (RS-2024-00345250)}
\thanks{\today}
\begin{abstract}
 In this paper we describe the geodesics on the K\"ahler cone of the Heisenberg group.  Furthermore we also prove that this is not a complete manifold.
\end{abstract}

\maketitle

\tableofcontents

\section{Introduction}

The Heisenberg group appears as one of the  most important objects in the study of not only many areas of mathematics but also in other fields such as physics and quantum mechanics. In particular, and for what concerns the purposes of this article, the one point compactification of the Heisenberg group is identified with the boundary of the complex hyperbolic plane and therefore it plays a vital role in the study of complex hyperbolic geometry.

The complex hyperbolic plane $\bH^2_\C$, endowed with the Bergman metric, is a K\"ahler manifold with constant holomorphic sectional curvature -1 and real sectional curvature pinched between -1 and -1/4. The group ${\rm SU}(2,1)$ is a triple cover of the group ${\rm PU}(2,1)$ of the holomorphic isometries of  $\bH^2_\C$. A full description of the geodesics of complex hyperbolic plane may be found in the book of Goldman, \cite{Go}, or, alternatively, in the notes of Parker, \cite{P}. 

Besides the Bergman metric, there is another natural metric which can be defined on the underlying manifold of complex hyperbolic plane. The full discussion about this metric may be found for instance in \cite{KPS}; we repeat here in brief its main features. The starting point of the construction of this metric is the Heisenberg group $\frakH$. The Heisenberg group is also a contact manifold and a CR-manifold. It turns out that from these structures it is defined a contact metric structure which is Sasakian. This is equivalent to say that the manifold $\calC(\frakH)=\R_{>0}\times\frakH$, that is, the horospherical parametrisation of complex hyperbolic plane, endowed with the warped product metric and a complex structure which is an extension of the almost complex structure of the horizontal space of $\frakH$ is a K\"ahler manifold. In contrast to the Bergman metric, it is of non constant negative holomorphic curvature. The manifold $\calC(\frakH)$ is the K\"ahler cone of the Heisenberg group and the goal of this paper is to describe its geodesics; incidentally, we do the same for the geodesics of the Heisenberg group as a Sasakian manifold. We also prove that the K\"ahler cone of the Heisenberg group is not complete.

The paper is organised as follows: in next section, we review the definition and properties of the Heisenberg group and its Riemannian cone constructed in \cite{KPS}. In section \ref{geo-hei}, we review the geodesics of $\frakH$ in \cite{Mar}. Section \ref{geo-cone} is intended to describe the geodesics on the Riemannian cone of the Heisenberg group. Finally,
we briefly prove that the Riemannian cone of the Heisenberg group is not geodesically complete.

\section{Preliminaries}\label{pre}
In this chapter, we briefly introduce the Heisenberg group and its properties. For more details we refer to \cite{KPS} as well as in \cite{Bo}, \cite{CDPT}, \cite{Go}, \cite{P}.

\subsection{Contact metric structure of the Heisenberg group}\label{sec-heis}
The (first) Heisenberg group $\frakH$  is the set $\C\times\R$ with multiplication $*$ given by
$$
(z,t)*(w,s)=(z+w,t+s+2\Im(z\overline{w})).
$$
The Heisenberg group $\frakH$ is a 2-step nilpotent Lie group; a left invariant basis for its Lie algebra comprises the following vector fields:
\begin{eqnarray*}
X=\frac{\partial}{\partial x}+2y\frac{\partial}{\partial t},\quad Y=\frac{\partial}{\partial y}-2x\frac{\partial}{\partial t},
\quad T=\frac{\partial}{\partial t}.
\end{eqnarray*}
The only nontrivial Lie bracket relation between them is
  $$
  [X,Y]=-4T.
  $$
A left-invariant 1-form $\omega$ is defined in $\frakH$ by 
\begin{equation*}\label{eq-omega}
\omega=dt+2xdy-2ydx=dt+2\Im(\overline{z}dz).
\end{equation*}
$\omega$ is a contact form for $\frakH$: $\omega\wedge d\omega\neq 0$ and in fact  $dm=-(1/4)\;\omega \wedge d\omega$, where $dm$ is the Haar measure for $\frakH$. 
The Reeb vector field for $\omega$ is $T$.
The kernel $\ker\omega=\calH$ is spanned by the vector fields $X,Y$.
Consider the almost complex structure $J$ defined on $\calH$ by
 $JX=Y$, $ JY=-X.$
 Then  $J$ is compatible with $d\omega$ and moreover, $\calH$ is a strictly pseudoconvex {\rm CR} structure; that is, $d\omega$ is positively oriented on $\calH$. 

 A sub-Riemannian metric $g_{cc}$ in $\frakH$ is defined in $\calH$ by the relations
$$
g_{cc}(X,X)=g_{cc}(Y,Y)=1,\quad g_{cc}(X,Y)=0.
$$
The corresponding sub-Riemannian tensor is then given by
$$
g_{cc}
=dx^2+dy^2.
$$
Given $L>0$ we are able to define a contact Riemannian structures in $\frakH$  such that the frame 
$$
\{X,Y, T_L=T/\sqrt{L}\}
$$ 
is orthonormal. In this way, we have that the Riemannian tensor is
$$
g_L
=dx^2+dy^2+L\omega^2=ds^2_{cc}+L\omega^2.
$$

Recall that if $(M,\omega)$ is a 3-dimensional contact manifold and $(\calH=\ker(\eta),J)$ is its CR-structure, then the almost complex structure $J$ on $\calH$ can be  extended to an endomorphism $\phi$ of the whole tangent bundle ${\rm T}(M)$ by setting $\phi(\xi)=0$, where $\xi$ is the Reeb field of $\omega$. Subsequently, a canonical Riemannian metric $g$ is defined in $M$ from the relations
\begin{equation*}\label{eq:contactmetric}
\omega(X)=g(X,\xi),\quad \frac{1}{2}d\omega(X,Y)=g(\phi X, Y),\quad \phi^2(X)=-X+\omega(X)\xi,
\end{equation*}
for all vector fields $X, Y$ in ${\mathfrak X}(M)$.
We then call $(M;\omega,\xi,\phi,g)$ the contact Riemannian structure  on $M$ associated to the contact structure $(M, \omega)$.

The complex structure $J$ on $\frakH$ can be extended to an endomorphism $\phi$ of the whole tangent bundle by setting $\phi(T_L)=0$.
Then it is proved in \cite{KPS} that the group $\frakH$ with contact form  $\sqrt{L}\omega$, Reeb vector field $T_L,$ endomorphism $ \phi$ and Riemannian metric $g_L$ is a contact Riemannian manifold if and only if $L=1/4$. 
Being in accordance with the notation in \cite{KPS}, we will hereafter write $g$ instead of $g_{1/4}$, $\widetilde{T}$ instead of $2T$ and $\widetilde{\omega}$ instead of $(1/2)\omega$.

\subsection{K\"ahler cone}\label{sec-Kahler}
With the contact metric structure described above, it is proved in \cite{KPS} that $\frakH$ is a Sasakian manifold. That is equivalent to say that the manifold
$\calC(\frakH)=\frakH\times\R_{>0}$ endowed with the metric
$$
g_r=dr^2+r^2g
$$
is K\"ahler. Explicitely:
\begin{itemize}
    \item An orthonormal basis for the metric $g_r$ comprises the vector fields
$$
X_r=(1/r)X,\quad Y_r=(1/r)Y,\quad T_r=(1/r)\widetilde{T},\quad R_r=d/dr.
$$
Note that all Lie brackets vanish besides
$$
[X_r,Y_r]=-(2/r)T_r,\quad[X_r,R_r]=(1/r)X_r,
\quad[Y_r,R_r]=(1/r)Y_r,\quad[T_r,R_r]=(1/r)T_r.
$$
\item The endomorphism  $\J$ defined by
$$
\J X_r=Y_r,\quad \J Y_r=-X_r,\quad \J T_r=-R_r,\quad \J R_r=T_r,
$$
is a complex structure for $\calC(\frakH)$ which preserves the metric $g_r$.
\item The 2-form $\Omega_r$ defined by
$$
\Omega_r=d\left(\frac{r^2}{2}\widetilde\omega\right)
$$
is the fundamental form for the Hermitian manifold $(\calC(\frakH), g_r, \J)$ and it is closed. Thus the manifold $(\calC(\frakH), g_r, \J, \Omega_r)$ is K\"ahler.
\end{itemize}

\section{Geodesics of the Heisenberg group}\label{geo-hei}
An exhaustive treatment of the geodesics of $(\frakH, g)$ may be found for instance in \cite{HZ}, \cite{Mar}, \cite{N}. We repeat in brief the discussion below. Let $\gamma(s)=(x(s),y(s),t(s))$ be a smooth curve defined in an interval $I=(-\epsilon,\epsilon)$, where $\epsilon>0$ and is sufficiently small, and suppose that $\gamma(0)=(x_0,y_0,t_0)=p_0$. 
The tangent vector $\dot \gamma(s)$ at a point $\gamma(s)$ is then
\begin{eqnarray*}
\dot\gamma(s)=\dot\gamma&=&\dot x\partial_x+\dot y\partial_y+\dot t\partial_t\\
&=& \dot xX+ \dot yY+(1/2)(\dot t+2x\dot y-2y\dot x)\widetilde{T}.
\end{eqnarray*}
We set
$$
f(s)=\dot x(s),\quad g(s)= \dot y(s),\quad h(s)= (1/2)(\dot t(s)+2x(s)\dot y(s)-2y(s)\dot x(s)).
$$
We may assume that $\gamma$ is of unit speed: $f^2+g^2+h^2=1$. 
Recall from \cite{KPS} that if $\nabla$ is the Riemannian connection then the following hold:
\begin{eqnarray*}
&&
\nabla_{X}X=0,\quad \nabla_{X}Y=-\widetilde{T},\quad \nabla_{X}\widetilde{T}=Y,\\
&&
\nabla_{Y}X=\widetilde{T},\quad\nabla_{Y}Y=0,\quad \nabla_{Y}\widetilde{T}=-X,\\
&&
\nabla_{\widetilde{T}}X=Y,\quad \nabla_{\widetilde{T}}Y=-X,\quad \nabla_{\widetilde{T}}\widetilde{T}=0.
\end{eqnarray*}
Recall again from \cite{KPS} that 
the covariant derivative of $\dot\gamma$ is
$$
\frac{D\dot\gamma}{ds}=\dot fX+\dot gY+\dot h\widetilde{T}+f\nabla_{\dot\gamma}X+g\nabla_{\dot\gamma}Y+
h\nabla_{\dot\gamma}\widetilde{T}.
$$
Since
\begin{eqnarray*}
&&
\nabla_{\dot\gamma}X=f\nabla_{X}X+g\nabla_{Y}X+h\nabla_{\widetilde{T}}X
=g\widetilde{T}+hY,\\
&&
\nabla_{\dot\gamma}Y=f\nabla_{X}Y+g\nabla_{Y}Y+h\nabla_{\widetilde{T}}Y=-f\widetilde{T}-hX,\\
&&
\nabla_{\dot\gamma}{\widetilde{T}}=f\nabla_{X}\widetilde{T}+g\nabla_{Y}\widetilde{T}+h\nabla_{\widetilde{T}}\widetilde{T}=fY-gX,
\end{eqnarray*}
we deduce
$$
\frac{D\dot\gamma}{ds}=(\dot f-2gh)X+(\dot g+2fh)Y+\dot h\widetilde{T}.
$$
Therefore, the geodesic equations are
\begin{equation}\label{eq:sys}
\dot f=2gh,\quad \dot g=-2fh,\quad \dot h=0,\quad f^2+g^2+h^2=1.
\end{equation}
In the special case $h=0$, that is, $\gamma$ is horizontal, we obtain the straight lines
\begin{equation*}\label{eq:geods}
\gamma(s)=(as+x_0,\;bs+y_0,\;2(ay_0-bx_0)s+t_0),
\end{equation*}
where $a,b$ are real constants and $a^2+b^2=1$. Those are all $g_{cc}$-geodesics. 

In the case $h\neq 0,$ one can get $h=c,$ where $c\neq0$ is a constant. We now write $F=f+ig,$ $z(s)=x(s)+iy(s)$ and $z_0=x_0+iy_0$. Then the above system (\ref{eq:sys}) is written equivalently as
$$
\dot F=-2ic\;F,
$$
and
has general solution
$$
F(s)=ke^{-2ics},\quad k\in\C, \quad |k|^2+c^2=1.
$$
We therefore have $|c|\le 1$. If $|c|=1$, then $k=0$ and
\begin{equation*}\label{eq:geodv}
\gamma(s)=(x_0,\;y_0,\;2cs+t_0),
\end{equation*}
is a vertical geodesics through $p_0$.
If now $|c|<1$ and $\gamma(s)=(z(s),t(s))$, then 
\begin{eqnarray}\label{eq:geod3}
z(s)&=&\frac{ik(e^{-2ics}-1)}{2c}+z_0,\\\label{eq:geod4}
t(s)&=&\frac{1}{c}\left((1+c^2)s-\frac{(1-c^2)\sin(2cs)}{2c}-\Re(\overline{z_0}k(e^{-2ics}-1))\right)+t_0.
\end{eqnarray}
We observe that the curves in (\ref{eq:geod3}) are Euclidean circles, which are the same as the projection of non linear $g_{cc}$-geodesics on the complex plane. However, Eqs. (\ref{eq:geod3}) and (\ref{eq:geod4}) do not give any non linear $g_{cc}$-geodesic. In fact, in order this to happen we should have $\dot t(s)+2\Im(\overline{z(s)}\dot z(s)=0$, that is, the geodesic must be horizontal. 
This condition leads to $h=0,$ i.e.,
$c=0$ which is impossible.


\section{Geodesics of the Riemannian cone}\label{geo-cone}

Let $\gamma(s) = (x(s),y(s), t(s), r(s))$ be a smooth curve in the K\"ahler cone and suppose that $\gamma(0) =(x_0, y_0, t_0, r_0) = q_0$. Its tangent vector is
\begin{eqnarray*}
\dot\gamma(s) = \dot\gamma&=& \dot x\partial_x + \dot y\partial_y + \dot t\partial_t +\dot rR_r\\
&=& r\dot xX_r + r\dot yY_r + (r/2)(\dot t+2x\dot y - 2y\dot x )T_r + \dot rR_r.
\end{eqnarray*}
We set
\begin{eqnarray*}
  && 
  f(s) = r(s)\dot x(s),\quad g(s) = r(s)\dot y(s),\\
&&
 h(s) =(1/2)r(s)(\dot t(s)+2x(s)\dot y(s)-2y(s)\dot x(s) ),\quad k(s) =\dot r(s), 
\end{eqnarray*}
and we may suppose that $f^2+g^2+h^2+k^2=1$. Recall again from \cite{KPS} that  if  $\nabla^r$ is the Riemannian connection, then
\begin{eqnarray*}
&&
\nabla^r_{X_r}X_r=-(1/r)R_r,\quad \nabla^r_{Y_r}X_r=(1/r)T_r,\quad \nabla^r_{T_r}X_r=(1/r)Y_r,\quad\nabla^r_{R_r}X_r=0,\\
&&
\nabla^r_{X_r}Y_r=-(1/r)T_r,\quad \nabla^r_{Y_r}Y_r=-(1/r)R_r,\quad \nabla^r_{T_r}Y_r=-(1/r)X_r,\quad\nabla^r_{R_r}Y_r=0,\\
&&
\nabla^r_{X_r}T_r=(1/r)Y_r,\quad \nabla^r_{Y_r}T_r=-(1/r)X_r,\quad \nabla^r_{T_r}T_r=-(1/r)R_r,\quad\nabla^r_{R_r}T_r=0,\\
&&
\nabla^r_{X_r}R_r=(1/r)X_r,\quad\nabla^r_{Y_r}R_r=(1/r)Y_r,\quad\nabla^r_{T_r}R_r=(1/r)T_r,\quad\nabla^r_{R_r}R_r=0.
\end{eqnarray*}
The covariant derivative of $\dot\gamma$ is 
\begin{eqnarray*}
\frac{D\dot\gamma}{ds}&=&\dot fX_r+\dot gY_r+\dot hT_r+\dot kR_r\\
&&+f\nabla^r_{\dot\gamma}X_r+g\nabla^r_{\dot\gamma}Y_r+h\nabla^r_{\dot\gamma}T_r+k\nabla^r_{\dot\gamma}R_r\\
&=&\dot fX_r+\dot gY_r+\dot hT_r+\dot kR_r\\
&&+(f/r)(-fR_r+gT_r+hY_r)\\
&&+(g/r)(-fT_r-gR_r-hX_r)\\
&&+(h/r)(fY_r-gX_r-hR_r)\\
&&+(k/r)(fX_r+gY_r+hT_r).
\end{eqnarray*}
From the vanishing of the covariant derivative and the unit speed assumption we obtain the following geodesic equations:
\begin{eqnarray}
\label{eq-geo1}
\dot f&=&(1/r)(2gh-kf),\\
\label{eq-geo2}
\dot g&=&(1/r)(-2fh-kg),\\
\label{eq-geo3}
\dot h&=&(1/r)(-kh),\\
\label{eq-geo4}
\dot k&=&(1/r)(1-k^2).
\end{eqnarray}
Equation (\ref{eq-geo4}) also reads as
$$
r\ddot r+(\dot r)^2=1.
$$
The positive solutions to this ODE are of the form
$$
r(s)=\sqrt{(s+c_1)^2+c_2},\quad c_1, c_2\in\R,\; c_2\ge0. 
$$
From the initial condition $r(0)=r_0$, we also have $c_1^2+c_2=r_0^2$; thus 
\begin{equation*}\label{eq-r}
r(s)=\sqrt{s^2+2c_1s+r_0^2},\;c_1\in\R,\;r_0^2-c_1^2\ge 0.
\end{equation*}
From Equation (\ref{eq-geo3}) we obtain
$$
h(s)=\frac{c_3}{r(s)}=\frac{c_3}{\sqrt{s^2+2c_1s+r_0^2}}, \quad c_3\in\R.
$$
We now have 
$$
f^2+g^2=1-h^2-k^2=\frac{r_0^2-c_1^2-c_3^2}{s^2+2c_1s+r_0^2}\ge 0.
$$
On the other hand, from
$$
c_3/r(s)=(1/2)r(s)(\dot t(s)+2x(s)\dot y(s)-2y(s)\dot x(s))
$$
we obtain that
\begin{equation}\label{eq-t0}
\dot t(s)+2x(s)\dot y(s)-2y(s)\dot x(s)=2c_3/(s^2+2c_1s+r_0^2).
\end{equation}

In the case where $r_0^2=c_1^2+c_3^2,$ one can get that $f\equiv g\equiv 0$ which yields to
$$
x(s)=x_0,\quad y(s)=y_0.
$$
Also, from $\dot t(s)=2\sqrt{r_0^2-c_1^2}/(s^2+2c_1s+r_0^2)$ we have
$$
t(s)=2\arctan\left(\frac{s\sqrt{r_0^2-c_1^2}}{r_0^2+c_1s}\right)+t_0,\quad r_0^2-c_1^2>0.
$$
In the case where $r_0^2=c_1^2$, it is easy to know that $r(s)=\pm s+r_0$ and $t(s)=t_0$. 
Hence the resulting geodesics in this case are of the form:
\begin{equation*}\label{eq-geoext}
\gamma_c(s)=\left(x_0,\;y_0,\;2\arctan\left(\frac{s\sqrt{r_0^2-c_1^2}}{r_0^2+c_1s}\right)+t_0,\sqrt{s^2+2cs+r_0^2}\right), \quad c\in\R
\end{equation*}
or straight lines of the form
\begin{equation*}\label{eq-geoextlines}
\gamma(s)=(x_0,\;y_0,\;t_0,\;s+r_0).
\end{equation*}
In the case where $r_0^2>c_1^2+c_3^2$, by plugging $r(s),h(s)$ and $k(s)$ into Eqs. (\ref{eq-geo1}) and (\ref{eq-geo2}) we obtain
\begin{eqnarray*}
\dot f&=&\frac{2c_3g-(s+c_1)f}{s^2+2c_1s+r_0^2},\\
\dot g&=&-\frac{2c_3f+(s+c_1)g}{s^2+2c_1s+r_0^2}.
\end{eqnarray*}
We set $F=f+ig$ and the system of geodesic equations becomes the following complex ODE of the first order:
\begin{equation*}
\dot F=-\frac{s+c_1+2ic_3}{s^2+2c_1s+r_0^2} F,
\end{equation*}
where
$$
|F|^2=\frac{r_0^2-c_1^2-c_3^2}{s^2+2c_1s+r_0^2}>0.
$$
Then we get
$$
F(s)=\frac{Ce^{i\Phi(s)}}{r(s)},
$$
where $C \in \mathbb{C}$ satisfies $|C|^2=r_0^2-c_1^2-c_3^2$ and 
\begin{equation}\label{eq-Phi}
   \Phi(s)=-\frac{2c_3}{\sqrt{r_0^2-c_1^2}}\arctan\left(\frac{s+c_1}{\sqrt{r_0^2-c_1^2}}\right).
\end{equation}

Since $F(s)=r(s)\dot z(s)$, we have the complex ODE of the first order
$$
\dot z(s)=\frac{Ce^{i\Phi(s)}}{r^2(s)},
$$
which gives
$$
z(s)=\frac{iCe^{i\Phi(s)}}{2c_3}+D,\quad D\in\C.
$$
From the initial conditions we then obtain
$$
D=z_0-\frac{iCe^{i\Phi(0)}}{2c_3}.
$$
Hence, we have
\begin{equation*}\label{eq-z}
z(s)=\frac{iC\left(e^{i\Phi(s)}-e^{i\Phi(0)}\right)}{2c_3}+z_0.
\end{equation*}
Now from  Eq. (\ref{eq-t0}) we have
\begin{eqnarray*}
    \dot t(s)&=&2\left(\frac{c_3}{r^2(s)}-\Im(\overline{z}(s)\dot z(s))\right)\\
&=&\frac{2c_3^2+|C|^2}{c_3r^2(s)}-\frac{|C|^2}{c_3r^2(s)}\cos\left(\Phi(s)-\Phi(0)\right)-\frac{2}{r^2(s)}\Im\left(\overline{z_0}C e^{i\Phi(s)}\right).
\end{eqnarray*}
By integrating and taking under account the initial conditions we obtain
\begin{eqnarray*}
t(s)&=&-\left(1+\frac{|C|^2}{2c_3^2}\right)(\Phi(s)-\Phi(0))+\frac{|C|^2}{2c_3^2}\sin(\Phi(s)-\Phi(0))\\
&&-\Re\left(\overline{z_0}C\left(e^{i\Phi(s)}-e^{i\Phi(0)}\right)\right)+t_0,\\
\end{eqnarray*}
where $\Phi(s)$ is given by (\ref{eq-Phi}).

\subsection{Non-completeness}\label{non}
In \cite{P-S}, the authors considered the non-completeness of the geodesics under the metric $g_{\calU}=dr^2+r^2dt^2$ on $\calU$, where $\calU=\{(t,r):t\in\R,\;r>0\}$ is the half plane.
It is also proved in \cite{KPS} that $\calU$  can be embedded into $C(\frakH)$ as a submanifold by setting 
$$\iota_{\calU}(t,r)=(0,2t,r )
,$$
which is totally geodesic with $g_\calU=\iota_{\calU}^{\ast}g_r.$ 
We now  obtain the following:
\begin{cor}
The totally geodesic submanifold  $\calU$ of $(\calC(\frakH),g_r)$ is not complete. Hence the same holds for $(\calC(\frakH),g_r)$.  
\end{cor}


\begin{thebibliography}{10}

  

\bibitem{Bo}
C.P.~Boyer,
\newblock {\em The {S}asakian geometry of the {H}eisenberg group}.
\newblock {Bull. Math. Soc. Sci. Math. Roumanie (N.S.)},
  52(100)(3):251--262, 2009.

  
  
  \bibitem{CDPT}
L.~Capogna and D.~Danielli and S.D.~Pauls and J.T. ~Tyson,
\newblock {\em An introduction to the {H}eisenberg group and the sub-{R}iemannian isoperimetric problem}.
\newblock Progress in Mathematics. Birkh\"{a}user Verlag, Basel,
  2007.
  
  

\bibitem{Go}
W.M.~Goldman,
\newblock {\em Complex hyperbolic geometry}.
\newblock Oxford University Press, 1999.

\bibitem{HZ}
P.~Hajłasz and S.~Zimmerman, 
\newblock {\em Geodesics in the Heisenberg Group}.
\newblock Anal. Geom. Metr. Spaces, 3(1): 325-337, 2015.

\bibitem{KPS}
J.~Kim and I.D.~Platis and L.~Sun,
\newblock  {\em PCR K\"ahler equivalent metrics in the Siegel domain}
\newblock Preprint, arXiv:2304.08079.




\bibitem{Mar}
V.~Marenich,
\newblock {\em Geodesics in Heisenberg groups}.
\newblock { Geom. Dedicata}, 66(2):175--185, 1997.

\bibitem{N}
G.A.~Noskov, 
\newblock{\em Geodesics in the Heisenberg group: an elementary approach}.
\newblock Sib. Elect. Math. Izv., 2008, Vol 5, 177–188.


\bibitem{P}
J.R.~Parker,
\newblock{\em Notes on Complex Hyperbolic Geometry}.
\newblock \texttt{https://maths.dur.ac.uk/users/j.r.parker/img/NCHG.pdf}


\bibitem{P-S}
I.D.~Platis and L.~Sun,
\newblock  {\em Half plane geometries: zero and unbounded negative curvature.}
\newblock Preprint, arXiv:2111.07569.








\end{thebibliography}
\end{document}